\documentclass[5p,times,authoryear]{elsarticle}
\pdfoutput=1
\usepackage{amsthm}
\usepackage{graphics}
\usepackage{amssymb,amsmath}
\usepackage[colorlinks,allcolors=blue]{hyperref}
\usepackage[utf8]{inputenc}
\usepackage[T1]{fontenc}
\usepackage{doi}
\usepackage{microtype}
\usepackage[capitalize]{cleveref}
\usepackage{hyphenat}
\usepackage{flushend}

\newcommand{\sdbst}{\ensuremath{\text{ST}_{\text{sd},b}}}
\newcommand{\poly}{\text{poly}}
\newcommand{\probname}[1]{\textsc{#1}}
\newcommand{\FPT}{\ensuremath{\text{FPT}}}
\newcommand{\APX}{\ensuremath{\text{APX}}}
\newcommand{\W}{\ensuremath{\text{W}}}
\newcommand{\N}{\ensuremath{\mathbb{N}}}
\newcommand{\pmtn}{\text{pmtn}}
\newcommand{\nwt}{\text{no-wait}}
\newcommand{\prc}{\text{prec}}
\newcommand{\chns}{\text{chains}}
\newcommand{\siz}{\text{size}}
\newcommand{\pP}{P}
\newcommand{\pR}{R}
\newcommand{\pO}{O}
\newcommand{\pJ}{J}
\newcommand{\pF}{F}
\newcommand{\NP}{\ensuremath{\text{NP}}}
\newcommand{\PP}{\ensuremath{\text{P}}}
\newcommand{\XP}{\ensuremath{\text{XP}}}

\newtheorem{openproblem}{Open Problem}

\begin{document}
\begin{frontmatter}
  \title{Parameterized complexity of machine scheduling: 15 open problems}
\author[bonn,mastricht]{Matthias Mnich}
\ead{mmnich@uni-bonn.de}

\author[nsu,imsoran]{René van Bevern\corref{cori}}
\ead{rvb@nsu.ru}

\address[bonn]{Universit{\"a}t Bonn, Bonn, Germany}
\address[mastricht]{Maastricht University, Maastricht, The Netherlands}
\address[nsu]{Department of Mechanics and Mathematics, Novosibirsk State University, Novosibirsk, Russian Federation}
\address[imsoran]{Sobolev Institute of Mathematics of the Siberian Branch of the Russian Academy of Sciences, Novosibirsk, Russian Federation}

\cortext[cori]{Corresponding author.  René van Bevern, Department of Mechanics and Mathematics, Novosibirsk State University, ul.\ Pirogova 2, 630090 Novosibirsk, Russian Federation}

\begin{abstract}
  \noindent
  Machine scheduling problems are a long-time key domain of algorithms and complexity research.
  A novel approach to machine scheduling problems are fixed\hyp parameter algorithms.
  To stimulate this thriving research direction,
  we propose %
  15 open questions in this area
  whose resolution we expect to lead
  to the discovery of new approaches and techniques
  both in scheduling and parameterized complexity theory.
\end{abstract}


\end{frontmatter}
\thispagestyle{empty}

\section{Introduction}
\label{sec:introduction}

\noindent
Algorithms for machine scheduling problems form one of the core applications of combinatorial optimization.
In those problems,
we are generally given a finite set~$J$ of jobs
with certain characteristics,
and we must find a \emph{schedule}
for processing the jobs on
one or more machines,
which also may have their individual specifications.
Typical characteristics of a job~$j$ are its processing time~$p_j\in\N$, its release date~$r_j\in \N$, its due date~$d_j\in\N$, or its importance reflected by a weight~$w_j\in\N$.
Jobs may be subject to precedence constraints
enforcing some jobs to be completed
before other jobs start.
Also, jobs may be preempted,
or may be required to be processed without preemption.
Machine characteristics typically include their speed or whether they are capable of processing a certain type of job.

Usually, one is not only
searching for a feasible schedule
that respects all constraints,
but additionally optimizes some objective function.
Classical objectives include the minimization of the makespan
or the sum of weighted completion times
(which is equivalent to minimizing the weighted average completion time).
Since the inception of the field in the 1950s,
thousands of research papers have been devoted
to understanding the complexity of scheduling problems.
A significant portion of investigated problems turned out to be \NP{}-hard.
In consequence, algorithm designers proposed algorithms for these problems that yield approximate solutions in polynomial time.
In 1999, \citet{SchuurmanWoeginger1999} listed 10 of the most prominent open problems around polynomial\hyp time approximation algorithms for \NP{}-hard scheduling problems at that time.

Only recently, a different algorithmic approach
has been put forward for solving \NP{}-hard scheduling problems:
\emph{fixed\hyp parameter algorithms}.
The idea in fixed\hyp parameter algorithms is
to accept exponential running times,
which are seemingly inevitable in solving \NP-hard problems,
but to restrict them to certain aspects of the problem,
which are captured by \emph{parameters} \citep{CyganEtAl2015,FlumGrohe2006,Niedermeier2006}.
More formally,
fixed\hyp parameter algorithms solve
an instance of size~\(n\) with parameter~\(k\)
in $f(k)\cdot \poly(n)$~time for some computable,
typically superpolynomial, function~$f$.
Thus,
fixed\hyp parameter algorithms
can solve even large instances
of \NP-hard scheduling problems
if the parameter takes only small values, the function $f$ grows only moderately,
and the polynomial degree in $n$ is small  \citep{vanBevernEtAl2015}.
Moreover,
problems that are even difficult to approximate
can be approximated efficiently and well
in real-world instances
using
fixed\hyp parameter \emph{approximation} algorithms
that exploit small parameters of real-world data
\citep{vanBevernEtAl2017}.

While fixed\hyp parameter algorithms are now a well\hyp investigated area of algorithmics,
their systematic application to scheduling problems
has gained momentum only recently
\citep[and more references below]{vanBevernPyatkin2016,ChenEtAl2017,CielibakEtAl2004,HalldorssonKarlsson2006,HermelinEtAl2015,HermelinEtAl2017b,JansenEtAl2017,MarxSchlotter2011}.
This already led to the transfer
of proof techniques
from parameterized complexity
to the world of scheduling,
such as
\(n\)-fold integer programming \citep{KnopKoutecky2017},
color coding,
problem kernelization \citep{vanBevernEtAl2015},
and
\W-hardness \citep{BentertEtAl2017,BodlaenderFellows1995,vanBevernEtAl2016,vanBevernEtAl2016b,MnichWiese2015,FellowsMcCartin2003,HermelinEtAl2017}.
It also led to the transfer of techniques from
mathematical programming to parameterized complexity,
such as convex integer programming
\citep{MnichWiese2015}
or parameterizing by structural properties
of the integer feasible polytope of linear programs
\citep{JansenKlein2017}.

In the following,
we summarize known results
and list open problems regarding the parameterized complexity of scheduling problems
on a single machine, on parallel identical machines,
and in shop scheduling environments.
We do not claim these problems to be the most important ones,
but expect that their resolution
(in one way or the other)
will lead to the discovery of further new approaches both in
parameterized complexity and scheduling theory,
thus stimulating further research
with both practical and theoretical significance.

In this sense,
we hope that this work will be appealing and inspiring both
to researchers with a scheduling background,
as well as to researchers with a parameterized complexity background.
For the latter,
we  exhibit some fixed\hyp parameter tractability results
that appeared before the advent of parameterized complexity
and thus were not explicitly described as such.

\paragraph{Organization of this work}
We start by giving preliminaries about machine scheduling and parameterized complexity in \cref{sec:preliminaries}.
Then we look at single\hyp machine problems in \cref{sec:singlemachineproblems}, and parallel identical machines in \cref{sec:parallelmachinescheduling}.
We finally consider the broad class of shop scheduling problems in \cref{sec:shopscheduling}.

\section{Preliminaries}
\label{sec:preliminaries}
\noindent
This section defines
basic notions of
scheduling theory,
parameterized complexity theory,
and approximation algorithms.

\subsection{Scheduling theory}
\noindent
Throughout this work, we use the standard three-field notation of scheduling problems due to \citet{GrahamEtAl1979}.
This allows us to denote many problems as a triple \(\alpha|\beta|\gamma\),
where \(\alpha\)~is the machine environment,
\(\beta\)~are job characteristics and scheduling constraints,
and \(\gamma\)~is the objective function.

\subsubsection{Machine models}
\noindent
We consider the following machine environments~\(\alpha\),
which are
described
in detail
in the following. %
Single\hyp machine environments are
denoted by the symbol~$1$,
parallel identical machines are denoted by~$P$,
job shop, open shop, and flow shop environments
are denoted by $J$, $O$, and $F$, respectively.

The symbol describing the machine environment
can be followed by an integer
restricting the number of machines
(for example, $P2$
indicates two parallel identical machines).
If the symbol is followed by $m$
(for example, $Pm$),
this indicates that the number of machines
is an arbitrary \emph{constant}.
If the symbol is neither followed
by a number nor~\(m\),
then the number of machines
is assumed to be given as part of the input.

\paragraph{Single-machine environments (denoted by ``1'' in the \(\alpha\)-field)}
Each job~\(j\) has to be processed for
a given amount~\(p_j\in\N\) of time
(its \emph{processing time})
on a single machine.
The machine can process only one job at a time.

\paragraph{Parallel identical machines (denoted by ``P'' in the \(\alpha\)-field)}
We are given a number~\(m\) of parallel identical machines
(that is, of the same speed).
Each job~\(j\) has to be processed by exactly one machine
for a given amount \(p_j\in\N\)~of time
and each machine can process only one job at a time.

\paragraph{Unrelated parallel machines (denoted by ``R'' in the \(\alpha\)-field)}
We are given a number~\(m\) of machines.
Each job~\(j\) has to be processed by exactly one machine~\(i\)
for a given amount \(p_{ij}\in\N\)~of time
and each machine can process only one job at a time.

\paragraph{Shop scheduling problems}
In shop scheduling problems,
we are given a set~$M$ of machines and
a set~$J$ of $n$~jobs,
where each job~$j$ consists of $n_j\in\N$~\emph{operations}.
In the most general setting,
the operations of each job are partially ordered:
an operation of a job can only start
once its preceding operations are completed.
Processing an operation~$o_{ij}$ of job~$j$,
where \(i\in\{1,\dots,n_j\}\),
requires time~$p_{ij}$
on a certain machine~$\mu_{ij}$.
Each job can be processed by
at most one machine at a time and
each machine can process
at most one operation at a time.
The three most important classes of shop scheduling problems are:

\emph{Open shop scheduling},
denoted by ``O'' in the \(\alpha\)-field:
the processing order of the operations of each job is unrestricted
and must be decided by an algorithm.
The only constraint is that each job has
exactly one operation on each machine.

\emph{Job shop scheduling},
denoted by ``J'' in the \(\alpha\)-field:
the operations of each job are totally ordered,
yet
each job may have a distinct total order on its operations.
Herein, several operations may require processing on the same machine
and not every job may have operations on all machines.

\emph{Flow shop scheduling},
denoted by ``F'' in the \(\alpha\)-field:
all jobs have the same set of operations
with the same total order.

\subsubsection{Job characteristics}
\noindent
For job characteristics, we add qualifiers to the $\beta$-field.
In this survey,
we consider the following types of job characteristics.

\paragraph{Precedence constraints}
If jobs are restricted
by precedence constraints,
then we add the qualifier \emph{``prec''}
to the \(\beta\)-field to
indicate that jobs are only allowed to start
after their predecessors 
are completed.
In the case where
the partial order induced by the precedence constraints
is the disjoint union of total orders,
we write \emph{``chains''} instead of ``prec''.

\paragraph{Processing time restrictions}
\looseness=-1
If all jobs~\(j\) have unit processing time,
we add ``\(p_j=1\)'' to the \(\beta\)-field.
If all jobs~\(j\) have an equal processing time~\(p\),
we add ``\(p_j=p\)'' to the \(\beta\)-field.
In case where processing times are bounded
from above by some constant~\(c\),
we add ``\(p_j\leq c\)'' to the \(\beta\)-field.
If the processing times
are restricted to values in some set~\(P\subseteq\mathbb N\),
we add ``\(p_j\in P\)'' to the \(\beta\)-field.

\paragraph{Machine restrictions}
If each job~\(j\) can be processed only on a subset~\(M_j\) of machines,
we add ``\(M_j\)''~to the \(\beta\)-field.

\paragraph{Release dates}
Each job~\(j\) may have a \emph{release date}~\(r_j\),
before which it cannot be started.
In this case,
we add ``\(r_j\)'' to the \(\beta\)-field.

\paragraph{Preemption}
Per default,
we assume that jobs are not allowed to be preempted.
Otherwise,
we add ``\emph{pmtn}'' to the \(\beta\)-field.

\paragraph{Multiprocessor jobs}
When a job~\(j\) may
occupy size$_j$~machines during execution,
we add the qualifier ``\textit{size}$_j$''
to the \(\beta\)-field;
such problems are known as ``multiprocessor scheduling''.
If the number of required machines
is restricted to values in some set~\(S\subseteq\mathbb N\),
we add ``\(\text{size}_j\in S\)'' to the \(\beta\)-field.

\paragraph{Batch setup times}
When jobs are partitioned
into %
\emph{batches}
and a sequence\hyp dependent setup time~\(s_{pq}\)
is needed when switching from a job of batch~\(p\)
to a job of batch~\(q\),
then,
in accordance with \citet{AllahverdiEtAl2008},
we add ``\sdbst{}'' to the \(\beta\)-field.

\subsubsection{Objective functions}
\noindent
The objective functions~$\gamma$ that we aim to minimize are
\begin{itemize}
    \item the \emph{makespan}~$C_{\max}=\max\{C_j\mid j\in J\}$, where \(C_j\) is the completion time of job~\(j\) in a schedule, and
    \item the weighted sum of completion times, \(\sum_{j\in J}w_jC_j\), when each job~\(j\in J\) also has a weight~\(w_j\in\N\).
\end{itemize}
If each job~\(j\in J\) has a due date~\(d_j\in\N\),
we also minimize
\begin{itemize}
    \item the total weighted tardiness~\(\sum_{j\in J}w_jT_j\), where
    \[T_j=\max\{0,\allowbreak d_j-C_j\},\text{ and }\]
    \item the weighted number of tardy jobs, \(\sum_{j\in J}w_jU_j\),
    where \(U_j=0\) if job~\(j\) is finished by its due date and \(U_j=1\) otherwise.
    
    This is equivalent to maximizing \emph{throughput}---the weighted number of jobs that get finished by their due dates.
\end{itemize}
Generally, we will drop the ``\(j\in J\)'' subscript under the sum
and we will refer to the unit\hyp weight variants of the scheduling problems by
simply dropping the \(w_j\) from the objective functions.

\subsection{Parameterized complexity}

\noindent
An instance of a \emph{parameterized problem}~\(\Pi\subseteq\Sigma^*\times\mathbb N\) is a pair~$(x,k)$ consisting of the \emph{input~$x$} and the \emph{parameter~$k$}.
A~parameterized problem~$\Pi$ is \emph{fixed\hyp parameter tractable~(FPT)}
if there is an algorithm solving any instance of~$\Pi$ with parameter~\(k\) and size~$n$ in $f(k) \cdot \poly(n)$~time for some computable function~$f$.
Such an algorithm is called a \emph{fixed\hyp parameter algorithm}.

Note that fixed\hyp parameter tractability for a parameter~\(k\)
is a much stronger property than polynomial\hyp time solvability
for constant~\(k\):
a fixed\hyp parameter algorithm running in \(O(2^k\cdot n)\)~time
takes polynomial time even for \(k\in O(\log n)\).
In contrast, an algorithm running in time \(O(n^k)\)
is polynomial for constant~\(k\),
but is often impractical
even for small values of $k$.

The main goal of parameterized complexity
is determining which
parameterized problems
have fixed\hyp parameter algorithms.
These can efficiently solve \NP-hard problems
in applications where the parameter takes only small values and the degree of the polynomial in $n$ is small.
Notably,
a problem may be fixed\hyp parameter tractable
with respect to one parameter,
but might be not with respect to another.
How to find parameters
that are both
small in applications
and
lead to fixed\hyp parameter algorithms,
is a research branch on its own \citep{FellowsEtAl2013,KomusiewiczNiedermeier2012,Niedermeier2010}.

\paragraph{Problem kernelization}
Parameterized complexity enabled a mathematical formalization of polynomial\hyp time data reduction with provable performance guarantees: \emph{kernelization}.
It has been successfully applied
to obtain effective polynomial\hyp time data reduction algorithms
for many \NP-hard problems \citep{GuoNiedermeier2007,Kratsch2014}
and also led to techniques for proving lower bounds
on the effectivity
of polynomial\hyp time data reduction
\citep{MisraEtAl2011,BodlaenderEtAl2014}.

A \emph{kernelization algorithm} for a
parameterized problem~\(\Pi\subseteq\Sigma^*\times\mathbb N\)
reduces an instance~$(x,k)$
to an instance~$(x',k')$
in polynomial time
such that
the size of~\(x'\) and~\(k'\) depends only on~$k$
and such that \((x,k)\in\Pi\) if and only
if \((x',k')\in\Pi\).
The instance~\((x',k')\) is called a \emph{problem kernel}.

Note that it is the parameter that allows us to measure the
effectivity of polynomial\hyp time data reduction
since an absolute statement like
``the data reduction shrinks the input by at least one bit''
for an \NP-hard problem
would imply that we can solve \NP-hard problems in
polynomial time.

\paragraph{Parameterized intractability}
To show that a problem is presumably not
fixed\hyp parameter tractable,
there is a parameterized analogue
of \NP\hyp hardness.
The parameterized analogue of \PP\({}\subseteq{}\)\NP{} is a hierarchy of complexity classes
\FPT\({}\subseteq{}\)\W[1]\({}\subseteq{}\)\W[2]\({}\subseteq{}\dots\subseteq{}\)\XP,
where \XP{} is the class of parameterized problems
that are solvable in polynomial time for constant parameter values.
The working hypothesis is that all inclusions in this hierarchy are proper, paralleling the famous $\PP \not=\NP$ conjecture.
A parameterized problem~$\Pi$ is
\emph{(weakly) W[$t$]-hard} for some \(t\in\mathbb N\)
if any problem in W[$t$] has a parameterized reduction to~\(\Pi\):
a \emph{parameterized reduction} from a parameterized problem~$\Pi_1$
to a parameterized problem~$\Pi_2$ is an algorithm mapping an instance~$(x,k)$
to an instance~$(x',k')$
in time~$f(k)\cdot\poly(|x|)$
such that $k'\leq g(k)$ and
\((x,k)\in\Pi_1\iff(x',k')\in\Pi_2\),
where \(f\)~and~\(g\) are arbitrary computable functions.
By definition,
no W[$t$]-hard problem is fixed\hyp parameter tractable
unless FPT${}={}$W[$t$].

A parameterized problem is
\emph{strongly $\W[1]$-hard}
if it is $\W[1]$-hard
even if all numbers in the input
are encoded in unary.

\subsection{Approximation algorithms}
\noindent
Since we only consider minimization problems,
we introduce approximation terminology
only for minimization problems.

An \emph{\(\alpha\)\hyp approximation}
is a solution to an optimization problem
with cost at most \(\alpha\)~times the cost of
an optimal solution.
A~\emph{polynomial\hyp time approximation scheme~(PTAS)} is an
algorithm that computes a \((1+\varepsilon)\)\hyp approximation
in polynomial time for any constant~\(\varepsilon>0\).
An \emph{efficient polynomial\hyp time approximation scheme~(EPTAS)}
is an algorithm that computes a
\((1+\varepsilon)\)\hyp approximation
in time \(f(1/\varepsilon)\cdot \poly(n)\)
for some computable function~\(f\) and~\(\varepsilon >0\).
Finally, a \emph{fully polynomial\hyp time approximation scheme (FPTAS)}
is an algorithm that computes a \((1+\varepsilon)\)\hyp approximation
in time \(\poly(n,1/\varepsilon)\) for any \(\varepsilon>0\).

It is known
that any problem having an EPTAS is
fixed\hyp parameter tractable
parameterized by the cost of an optimal solution \citep[Lemma 11]{CesatiTrevisan1997}.
Thus,
showing $\W[1]$\hyp hardness of a problem
parameterized by the optimum solution cost
shows that it neither has an EPTAS nor an FPTAS
unless $\FPT{}=\W[1]$.
Moreover,
strongly \NP{}-hard optimization problems
with polynomially bounded objective functions
do not have FPTASes
unless $\PP=\NP$ \citep{GareyJohnson1979}.

\section{Single-machine problems}
\label{sec:singlemachineproblems}
\noindent
In this section,
we consider various single-machine problems:
\cref{sec:single:twt}
studies total weighted tardiness as objective function,
\cref{sec:single:twct} the total weighted completion time,
\cref{sec:single:wtj} the weighted number of tardy jobs,
and, finally,
\cref{sec:forbidden} considers a variant
with forbidden start and end times of jobs.

\subsection{Total weighted tardiness}
\label{sec:single:twt}
\noindent
The problem 1||\(\Sigma w_jT_j\)
is strongly NP-hard \citep{Lawler1977,LenstraEtAl1977},
yet has two easier special cases:
1|$p_j{=}p$|$\Sigma w_jT_j$
is solvable in polynomial time
via a reduction to the classical linear assignment problem,
whereas
1||\(\Sigma T_j\) is solvable
in pseudo\hyp polynomial time \citep{Lawler1977,LenstraEtAl1977}.

To be more exact,
\citet{Lawler1977}'s pseudo polynomial\hyp time algorithm
is a fixed\hyp parameter algorithm
for~1||\(\Sigma T_j\)
parameterized by the maximum processing time~\(p_{\max}\)
that also works for %
1||$\Sigma w_jT_j$
with \emph{agreeable} processing times and weights,
that is, \(p_i<p_j\) implies \(w_i\geq w_j\).

Although %
1||\(\Sigma w_jT_j\)
is strongly \NP-hard and thus
cannot have pseudo polynomial\hyp time algorithms
unless \PP${}={}$\NP{},
it is interesting whether
the fixed\hyp parameter tractability result
for  1||\(\Sigma w_jT_j\) with
agreeable processing times and weights holds in general:

\begin{openproblem}
  Is $1||\Sigma w_jT_j$ fixed\hyp parameter tractable parameterized by the maximum processing time~\(p_{\max}\)?
\end{openproblem}

\noindent
Also note that,
given that %
1|$p_j{=}p$|$\Sigma w_jT_j$
is solvable in polynomial time,
it is interesting to see
whether 1||$\Sigma w_jT_j$
is fixed\hyp parameter tractable with respect to
the number~\(\bar p\) of
distinct processing times.

\subsection{Total weighted completion time}
\label{sec:single:twct}
\noindent
\citet{AmbuhlMastrolilli2009} showed that
$1|\prc{}|\Sigma w_jC_j$ is a special case of \probname{Weighted Vertex Cover}.
\probname{Vertex Cover} is one of the most well\hyp studied
problems in parameterized complexity theory,
in particular in terms of
problem kernels; small kernels for \probname{Weighted Vertex Cover} were established by \citet{EtscheidEtAl2017}.
It is interesting
which
  kernelization algorithms
  carry over to $1|\prc{}|\Sigma w_jC_j$.
However,
a problem kernel for $1|\prc{}|\Sigma w_jC_j$ whose
size is bounded by a function of the optimum
is not interesting---the optimum is at least the weighted sum of all processing times
and thus already bounds the input size without data reduction.
Remarkably,
\probname{Vertex Cover} admits a (randomized) algorithm
yielding a problem kernel with size polynomial
in the
difference between the optimum
and a lower bound given by the relaxation of the natural ILP \citep{KratschWahlstrom2012}.
Since $1|\prc{}|\Sigma w_jC_j$ is a special case of \probname{Weighted Vertex Cover},
and thus allows for a likewise natural ILP formulation \citep{AmbuhlMastrolilli2009},
the following question is interesting:
\begin{openproblem}
  Does $1|\prc{}|\Sigma w_jC_j$ admit a problem kernel
  of size polynomial in the difference between the optimum
  and the lower bound obtained from its natural LP relaxation?
\end{openproblem}
\noindent
We point out that any partial progress on this question
would be interesting: be it a randomized kernelization algorithm
or even a partial kernel that bounds only the number of jobs
and not necessarily their weights or processing times.

The question is complicated by the fact
that each vertex in the \probname{Weighted Vertex Cover}
instance created by \citet{AmbuhlMastrolilli2009},
and thus each variable in the corresponding ILP,
corresponds to a \emph{pair} of jobs in the $1|\prc{}|\Sigma w_jC_j$ instance,
such that known data reduction rules for \probname{Weighted Vertex Cover}
do not allow for a straightforward interpretation in terms of jobs.

\subsection{Throughput or weighted number of tardy jobs}

\label{sec:single:wtj}
\label{sec:throughputsched}

\noindent
In a survey
on open questions in maximum throughput scheduling,
\citet{Sgall2012} asked
whether there is a polynomial\hyp time algorithm for
\(1|r_j,p_j{\leq}c|\Sigma U_j\)
for constant~\(c\).
Similarly,
the \NP{}\hyp hardness of the weighted
case for constant~\(c\) is open.
Parameterized complexity can serve
as an intermediate step towards resolving
these questions.
\begin{openproblem}\label{sgall}
  Are
  \(1|r_j|\Sigma U_j\) and \(1|r_j|\Sigma w_jU_j\)
  $\W[1]$-hard
  parameterized by the maximum processing time~$p_{\max}$?
\end{openproblem}
\noindent
Currently, even containment in the parameterized complexity class \XP{} is open.
The special case \(1||\Sigma U_j\) %
is polynomial\hyp time solvable,
whereas %
\(1||\Sigma w_jU_j\)
is weakly \NP{}-hard and solvable in pseudo\hyp polynomial time \citep{LawlerMoore1969,Karp1972}.

\citet{HermelinEtAl2017b} gave fixed\hyp parameter algorithms
for \(1||\Sigma w_jU_j\) simultaneously parameterized by any two
out of the following three parameters:
the number of distinct due dates, the number of distinct processing times, and the number of distinct job weights.

\citet{FellowsMcCartin2003}
showed %
that \(1|\prc{},p_{j}{=}1|\Sigma U_i\)
is \W[1]-hard parameterized by the number of tardy jobs
but fixed\hyp parameter tractable
with respect to this parameter
if the partial order induced by the precedence constraints
has constant width.
Herein, the width of a partial order is
the size of a largest set of mutually incomparable jobs.

\subsection{Forbidden start and end times}
\label{sec:forbidden}
\noindent
Machine scheduling problems with forbidden start and end times
model the situation when an additional resource,
subject to unavailability constraints,
is required to start or finish a job.
For example, this resource might be operators of chemical experiments, which serves as a major motivation for such problems.

For makespan minimization on a single machine,
\citet{BillautSourd2009} gave an algorithm
that runs in $n^{O(\tau^2)}$~time
for $\tau$~forbidden start times
and $n$~jobs;
this was improved by \citet{RapineBrauner2013} to $n^{O(\tau)}$~time.
For the high\hyp multiplicity encoding of the input---given by binary numbers~$n_t$ encoding the number of jobs having the same forbidden start and end times---\citet{GabayEtAl2016} showed a polynomial\hyp time algorithm if the number~$\tau$ of forbidden times is constant.
All of these results leave open the possibility for fixed\hyp parameter tractability of the problem parameterized by~\(\tau\).
\begin{openproblem}
  Is
  makespan minimization on a single machine
  with \(\tau\)~forbidden start and end times
  fixed\hyp parameter tractable parameterized by~$\tau$?
\end{openproblem}

\section{Parallel identical machines}
\label{sec:parallelmachinescheduling}
\noindent

\noindent
In this section,
we consider scheduling on parallel identical machines.
\cref{sec:parallelfewproctimes}
considers the makespan objective,
\cref{sec:prec-sched-fpt} the makespan objective
with precedence constraints,
\cref{sec:tardyjobs} the weighted number of tardy jobs,
\cref{sec:jit} considers just-in-time scheduling,
and, finally,
\cref{sec:preempt} considers variants
with allowed preemption of jobs.
\subsection{Makespan}
\label{sec:parallelfewproctimes}
\noindent
\citet{AlonEtAl1998} showed an EPTAS for $\pP{}||C_{\max}$,
which implies that $\pP{}||C_{\max}$ is
fixed\hyp parameter tractable parameterized by the makespan~$C_{\max}$
\citep[Lemma 11]{CesatiTrevisan1997}.

This %
was improved by \citet{MnichWiese2015}
to a fixed\hyp parameter algorithm for %
the maximum processing time~\(p_{\max}\) of any job;
we generally expect $p_{\max} \ll C_{\max}$.
The running time of both algorithms is doubly\hyp exponential in the parameter.
An improved algorithm whose running time depends singly\hyp exponentially on~$p_{\max}$ was proposed by \citet{KnopKoutecky2017}.
Very recently, \citet{KnopEtAl2017} strengthened this result by giving an algorithm with single\hyp exponential dependence $p_{\max}$ even for the \emph{high-multiplicity encoding}, when, for each processing time~$p\leq p_{\max}$, the number of jobs with processing time~$p$ is encoded in binary. 
\citet{ChenEtAl2017}
generalized the result of \citet{MnichWiese2015}
by showing that \(\pR{}||C_{\max}\)
is fixed\hyp parameter tractable parameterized by~\(p_{\max}\)
and the rank of the matrix~\((p_{ij})\)
giving the processing time~\(p_{ij}\) of job~\(j\)
on machine~\(i\).

In a breakthrough result,
\citet{GoemansRothvoss2014} %
showed that $\pP{}||C_{\max}$ is polynomial\hyp time solvable
when the number~$\overline{p}$ of distinct processing times
is constant, that is, the problem is in $\XP$.
Their result holds even
for the high-multiplicity encoding of the input,
where the number~$m$ of machines and
the number~$n_j$ of jobs with processing time~$p_j$
are encoded in binary
for each~\(j\in\{1,\dots,\overline p\}\).
A close inspection of their result reveals that it is a fixed-parameter algorithm when the processing times~$p_j$ are encoded in unary.
Despite all this progress,
the following problem remains open.
\begin{openproblem}
    Is $\pP{}||C_{\max}$ with high-multiplicity encoding fixed\hyp parameter tractable parameterized by the number~$\overline{p}$ of distinct processing times, that is, when job processing times~$p_j$ and the number of jobs for each processing time are encoded in binary?
\end{openproblem}
\noindent
To find a fixed\hyp parameter algorithm for the high-multiplicity encoding of $\pP{}||C_{\max}$ parameterized by~\(\overline{p}\), one difficulty is
outputting an optimal schedule,
whose obvious encoding requires at least \(m\)~bits
(to store how many jobs of each processing time are
processed on each machine)
and is therefore neither polynomial in the input size
nor bounded by a function of~\(\overline p\).
Such problems might be possibly overcome by using the framework of \citet{BraunerEtAl2005}.

\subsection{Makespan and precedence constraints}
\label{sec:prec-sched-fpt}
\noindent
The parameterized complexity of $\pP|\prc{}|C_{\max}$
and special cases has extensively been studied
with respect to the
width of the partial order
induced by the precedence constraints:
the \emph{width} of a partial order is the size
of its largest antichain---a set of pairwise incomparable jobs.

The special case
$\pP2|\chns{}|C_{\max}$
is weakly \NP{}-hard even
for partial orders of width three
\citep{GuntherEtAl2014,vanBevernEtAl2016}.
Since this excludes fixed\hyp parameter algorithms
using the partial order width as parameter
already on two machines,
it has been tried to use the partial order width
and the maximum processing time as parameters simultaneously.
\citet{BodlaenderFellows1995}
showed that
even %
$\pP|\prc{},\allowbreak p_j{=}1|C_{\max}$
is $\W[2]$-hard parameterized by the partial order width
and by the number of machines.
Later,
\citet{vanBevernEtAl2016}
showed that combining partial order width
with maximum processing time
does not even yield fixed\hyp parameter algorithms
on two machines.
More precisely,
they showed that even %
$\pP2|\prc{},\allowbreak p_j{\in}\{1,2\}|C_{\max}$
is $\W[2]$-hard parameterized by the partial order width;
and so is %
$\pP3|\prc{},\siz_j{\in}\{1,2\}|C_{\max}$.

Further restricting
this problem,
one
arrives at
a long\hyp standing open problem due to \citet[OPEN8]{GareyJohnson1979}
of whether %
$\pP3|\prc{},p_j {=} 1|C_{\max}$
is \NP{}-hard or polynomial\hyp time solvable.
In fact, the complexity is open for any constant number of machines.
Thus,
as pointed out by \citet{vanBevernEtAl2016},
it would be surprising to show $\W[1]$\hyp hardness of this problem for \emph{any} parameter,
since this would exclude polynomial\hyp time solvability unless $\FPT{} = \W[1]$.

\begin{openproblem}
  Is $\pP3|\prc{},p_j {=} 1|C_{\max}$ fixed\hyp parameter tractable parameterized by  the width of the partial order induced by the precedence constraints?
\end{openproblem}

\noindent
Note that a negative answer would basically
answer the open question of \citet{GareyJohnson1979}
on whether the problem is polynomial\hyp time solvable,
whereas a positive answer would
be in strong contrast to the $\W[2]$\hyp hardness
of the slight generalizations
$\pP|\prc{},\allowbreak p_j{=}1|C_{\max}$,
$\pP2|\prc{},\allowbreak p_j{\in}\{1,2\}|C_{\max}$,
and
$\pP3|\prc{},\siz_j{\in}\{1,2\}|C_{\max}$
considered by \citet{BodlaenderFellows1995} and \citet{vanBevernEtAl2016}.

\subsection{Throughput or weighted number of tardy jobs}
\label{sec:tardyjobs}
\noindent
\noindent
\citet{BaptisteEtAl2004}
showed that
\(\pP m|r_j,p_j{=}p|\Sigma w_jU_j\)
is polynomial\hyp time solvable.
According to \citet{Sgall2012},
this is open when
the number of machines is not a constant.
A first step to resolving this question
is solving the following problem.
\begin{openproblem}
  Are \(\pP|r_j,p_j{=}p|\Sigma w_jU_j\)
  and \(\pP|r_j,p_j{=}p|\Sigma U_j\)
  fixed\hyp parameter tractable parameterized by the number of machines?
\end{openproblem}
\noindent
A negative answer to this question
would also be interesting
since it is open whether these problems are even \NP{}-hard
if the number~\(m\) of machines is part of the input.
Notably,
the special case \(\pP|p_j{=}p|\Sigma w_jU_j\)
is polynomial\hyp time solvable
via a reduction to the classical linear assignment problem
but %
\(\pP|p_j{=}p,\pmtn|\Sigma w_jU_j\)
is strongly \NP{}-hard \citep{BruckerKravchenko1999}.

\subsection{Interval scheduling or just-in-time scheduling}
\label{sec:jit}
\noindent
An important special case
of maximizing the throughput on parallel identical machines
is \emph{interval scheduling},
where each job~\(j\) has a weight~\(w_j\),
a fixed start time, and a fixed end time.
The goal is to schedule a maximum\hyp weight subset of jobs
non\hyp preemptively on \(m\)~parallel identical machines.
As always, each machine can process only one job at a time.
Since we can interpret this problem as minimizing the total weight of jobs not meeting their due dates, where each job~\(j\) has a processing time~\(p_j\), a release date~\(r_j\), and a due date~\(d_j\) such that \(d_j-r_j=p_j\), we denote the problem as \(\pP|d_j{-}r_j{=}p_j|\Sigma w_jU_j\).
This problem is
equivalent to maximizing
the weighted number
of ``just-in-time'' jobs
on parallel identical machines
\citep{SungVlach2005}.

\citet{ArkinSilverberg1987}
show that this problem is solvable in $O(n^2\log n)$~time.
However,
they showed that
the variant \(\pP|M_j,d_j{-}r_j{=}p_j|\Sigma w_jU_j\)
is \NP{}-hard and
solvable in $O(n^{m+1})$~time.
\begin{openproblem}\label{cis}
  Is %
  \(\pP|M_j,d_j{-}r_j{=}p_j|\Sigma w_jU_j\)
  fixed\hyp parameter tractable parameterized by the number of machines?
\end{openproblem}
\noindent
This problem
is a special case
of \(\pR|d_j{-}r_j{=}p_j|\Sigma w_jU_j\)
where each job~$j$
has a processing time~$p_{ij}\in\mathbb \{p_j,\infty\}$ on machine~\(i\).
In that sense, any positive answer to \cref{cis} (in form of a fixed\hyp parameter algorithm) can serve as a first step towards
a fixed\hyp parameter algorithm for
\(\pR|d_j{-}r_j{=}p_j|\Sigma w_jU_j\) also;
after all,
the \(O(n^{m+1})\)-time algorithm of \citet{ArkinSilverberg1987}
for   \(\pP|M_j,d_j{-}r_j{=}p_j|\Sigma w_jU_j\)
was generalized to
\(\pR|d_j{-}r_j{=}p_j|\Sigma w_jU_j\)
with essentially the same running time \citep{SungVlach2005}.

We point out that, if we only slightly relax
the condition \(d_j-r_j=p_j\)
to \(d_j-r_j\leq\lambda p_j\) for any constant~\(\lambda>0\),
then the problem is weakly \NP{}-hard for \(m=2\) and
strongly \W[1]-hard parameterized by~\(m\) \citep{vanBevernEtAl2016b},
even when all jobs are allowed to be processed on all machines
and
only checking whether one can finish \emph{all}
jobs by their due date.

Using a construction due to \citet{HalldorssonKarlsson2006},
\(\pP|M_j,d_j{-}r_j{=}p_j|\Sigma w_jU_j\)
can be seen to be a special case of the \probname{Job
  Interval Selection} problem introduced by \citet{NakajimaHakimi1982}, where each job has multiple
possible execution intervals
and we have to select one execution interval for each job
in order to process it: we model the machines by pairwise disjoint segments of the real line and each job has an interval in each segment
belonging to a machine it can be processed on.
Via this relation to \probname{Job Interval Selection},
one can model \(\pP|M_j,d_j{-}r_j{=}p_j|\Sigma w_jU_j\)
as \probname{Colorful Independent Set} 
in an interval graph 
with at most \(mn\)~intervals colored in \(n\)~colors and having at most \(m\)~intervals of each color
\citep{vanBevernEtAl2015},
where the task is to find a maximum\hyp weight independent
set of intervals with mutually distinct colors.

\probname{Colorful Independent Set}
is fixed\hyp parameter tractable parameterized
by the number of colors in the solution
on interval (and even on chordal)
graphs \citep{vanBevernEtAl2015,BentertEtAl2017},
which implies
that \(\pP|M_j,d_j{-}r_j{=}p_j|\Sigma w_jU_j\) is fixed\hyp parameter tractable
parameterized by the number of jobs that can be scheduled in an optimal solution.
To answer \cref{cis},
one could try to show that
\probname{Colorful Independent Set} in interval graphs
is fixed\hyp parameter tractable parameterized by the maximum number of intervals of any color.

\subsection{Allowed preemption}
\label{sec:preempt}
\noindent
The three problems
\(\pP|\pmtn{},r_j|\Sigma C_j\),
\(\pP|\pmtn|\Sigma U_j\),
and
\(\pP|\pmtn|\Sigma T_j\)
are all NP-hard
\citep{DuEtAl1990,Lawler1983,KravchenkoWerner2013},
yet the special cases with equal processing times
are polynomial\hyp time solvable
\citep{BaptisteEtAl2004,BaptisteEtAl2007,KravchenkoWerner2011}.

\begin{openproblem}
  Are the above problems fixed\hyp parameter tractable parameterized by the number $\overline{p}$ of distinct processing times?
\end{openproblem}
\noindent
\citet{KravchenkoWerner2011} survey
related problems;
as of now, the complexity status of %
$\pP2|\pmtn{},p_j {=} p|\Sigma w_jT_j$ 
for constant~\(p\) is open.
Parameterized hardness might serve as an intermediate step
towards settling it. %
\begin{openproblem}
    Is $\pP2|\pmtn{},p_j {=} p|\Sigma w_jT_j$ \W[1]-hard parameterized by the processing time~$p$ of all jobs?
\end{openproblem}

\section{Shop scheduling}
\label{sec:shopscheduling}
\noindent

\noindent
We now consider the open shop, job shop, and flow shop scheduling problems.
\cref{sec:jcmax} considers the makespan objective,
\cref{sec:cmaxapr} considers fixed\hyp parameter approximation schemes,
\cref{sec:setup} takes into account setup times,
and, finally,
\cref{sec:jsthroughput} considers variants of maximizing throughput.

\subsection{Makespan}\label{sec:jcmax}
\noindent
\citet{KononovEtAl2012} give a computational complexity dichotomy
of \(\pJ||C_{\max}\) and \(\pO||C_{\max}\)
into
polynomial\hyp time solvable cases and \NP{}-hard cases
depending on the maximum processing time \(p_{\max}\) of operations,
the maximum number~\(n_{\max}\) of operations per job,
and an upper bound on the makespan~\(C_{\max}\).
All problems are \NP{}-hard even if
all of the listed parameters are simultaneously
bounded from above by~4 \citep{KononovEtAl2012,WilliamsonEtAl1997},
which fully settles the parameterized complexity
of \(\pO||C_{\max}\) and \(\pJ||C_{\max}\)
with respect to these parameters.

However,
in all hardness reductions of \citet{KononovEtAl2012},
the number of machines is necessarily unbounded:
when bounding both the number~\(m\) of machines and the makespan~\(C_{\max}\),
then
shop scheduling problems are trivial,
since the overall number of operations
in the input is at most~\(m\cdot C_{\max}\).
This makes the parameters~\(p_{\max}\) and \(n_{\max}\)
interesting for fixed\hyp parameter algorithms
for shop scheduling problems
with a fixed number of machines
(or with the number of machines as an additional parameter).

For example,
both %
\(\pF 3||C_{\max}\) and
\(\pO 3||C_{\max}\)
are strongly \NP{}-hard \citep{GonzalezSahni1976,GareyEtAl1976},
yet
\(\pO m||C_{\max}\) is fixed\hyp parameter
tractable parameterized by
maximum processing time~\(p_{\max}\):
if the maximum machine load is
\(\ell_{\max}\geq m^2 p_{\max}\),
then the makespan of the instance is \(\ell_{\max}\)
\citep{Sev95}.
Otherwise,
there are at most
\(\ell_{\max}< m^2 p_{\max}\)~jobs
and the instance can be solved using brute
force.
For \(\pF 3||C_{max}\),
the analogous question is open:

\begin{openproblem}\label{prob:o3}
  Is \(\pF3||C_{\max}\) fixed\hyp parameter tractable
  parameterized by the maximum processing time~\(p_{\max}\)?
\end{openproblem}

\noindent
This question is likewise interesting for
\(\pF 3|\nwt|C_{\max}\),
where ``\nwt'' means that
each operation of a job has to start immediately
after completion of the preceding operation of the same job.
The \NP{}\hyp hardness of
\(\pF 3|\nwt|C_{\max}\)
was a long\hyp standing open question \citep{Rock1984}.

Notably,
for the three\hyp machine job shop scheduling problem,
the above question has a negative answer:
even the special case
\(\pJ3|p_{ij}{=}1|C_{\max}\) %
is \NP{}-hard \citep{LenstraRinnooyKan1979}.

\begin{openproblem}
  Is \(\pJ3||C_{\max}\) fixed\hyp parameter tractable
  parameterized by the maximum processing time~\(p_{\max}\)
  and the maximum number~\(n_{\max}\) of operations per job?
\end{openproblem}

\subsection{Makespan approximation}
\label{sec:cmaxapr}
\noindent
Job and open shop scheduling problems
are often \APX{}-hard \citep{WilliamsonEtAl1997}
but also \NP{}-hard for
a constant number of machines
\citep{GonzalezSahni1976,GareyEtAl1976,LenstraRinnooyKan1979}
and many other constant parameters \citep{KononovEtAl2012}.
Thus,
they are amenable neither to approximation schemes
nor fixed\hyp parameter algorithms.

However,
while having a constant number of machines does not help to get polynomial\hyp time algorithms
for these problems,
it helps tremendously in getting PTASes \citep{SevastianovWoeginger1998,Hall1998,JansenSolisOba2003}.
Remarkably,
these are not simply PTASes
for a fixed number of machines running in time,
say \(n^{O(m/\varepsilon)}\):
In terms of \citet{Marx2008},
they are actually
\emph{fixed\hyp parameter tractable approximation schemes} (\emph{FPT-AS}) 
for the parameters~\(m\) and~\(\varepsilon\),
running in time \(f(m,\varepsilon)\cdot\poly(n)\).
While not necessarily being a PTAS,
an FPT-AS
using the number~\(m\) of machines as parameter and
running in \(O(2^{m/\varepsilon}\cdot n)\)~time
might be practically more valuable
than a PTAS running in \(O(n^{1/\varepsilon})\)~time
or even an FPTAS running in \(O((n+1/\varepsilon)^3)\)~time.

\citet{SevastianovWoeginger1998} compute a \((1+\varepsilon)\)\hyp approximation for \(\pO||C_{\max}\) in \(f(m,\varepsilon)+O(n \log n)\)~time,
\citet{Hall1998} computes a \((1+\varepsilon)\)\hyp approximation for \(\pF||C_{\max}\) in \(f(m,\varepsilon)+O(n^{3.5})\)~time.
Yet for  \(\pJ||C_{\max}\),
\citet{JansenSolisOba2003}
need an additional parameter:
they compute a \((1+\varepsilon)\)\hyp approximation for \(\pJ||C_{\max}\)
in \(f(m,n_{\max},\varepsilon)+O(n)\)~time, where \(n_{\max}\)~is the maximum number of operations of a job.

\begin{openproblem}
  Is there an algorithm that yields a \((1+\varepsilon)\)\hyp approximation for \(\pJ||C_{\max}\)
  in \(f(m,\varepsilon)\cdot \poly(n)\)~time?
\end{openproblem}
\noindent
Giving a negative answer to this question seems challenging:
the obvious approach would be proving that \(\pJ||C_{\max}\)
is $\W[1]$-hard parameterized by~\(C_{\max}+m\).
However, as discussed in \cref{sec:jcmax},
\(\pJ||C_{\max}\) is trivially fixed\hyp parameter tractable parameterized
by~\(C_{\max}+m\),
so that this approach is inapplicable.

\subsection{Makespan with sequence\hyp dependent setup times}
\label{sec:setup}
\noindent
In the following, we consider
a variant \(\pO|p_{ij}=1,\sdbst|C_{\max}\)
of open shop
with unit processing times
where
the jobs are partitioned into \(g\)~batches
and machines require a
sequence\hyp dependent
batch setup time~\(s_{pq}\in\N\)
when switching from jobs of batch~\(p\)
to jobs of batch~\(q\).

The case with unit processing times
models applications where large batches of items
have to be processed and processing individual items
in each batch takes significantly less time than
setting up the machine for the batch.
Alternatively
one can interpret this problem as an open shop problem
where machines or experts have to travel between jobs in different locations
in order to process them \citep{AverbakhEtAl2006}.

As a generalization of the \probname{Traveling Salesperson} problem,
\(\pO|p_{ij}=1,\sdbst|C_{\max}\) is \NP{}-hard already for one machine.
As shown by \citet{vanBevernPyatkin2016},
the problem is solvable in \(O(2^gg^2+mn)\)~time
if each batch contains more jobs than there are machines.
For the case where batches may contain less jobs,
they showed that
the problem is solvable in \(O(n\log n)\)~time
if both the number~\(m\) of machines
and the number~\(g\) of batches are constants.
\begin{openproblem}
  Is \(\pO|p_{ij}=1,\sdbst|C_{\max}\) fixed\hyp parameter tractable parameterized by the number~\(g\) of batches?
\end{openproblem}

\noindent
We point out that, for this problem, it would be desirable to design a fixed\hyp parameter algorithm that works also for the high\hyp multiplicity encoding,
which compactly encodes
the number of jobs in each batch
in binary.
The fixed\hyp parameter algorithms of \citet{vanBevernPyatkin2016}
also work for the high\hyp multiplicity encoding as long as they only need to compute the minimum makespan.

\subsection{Throughput or weighted number of tardy jobs}
\label{sec:jsthroughput}
\noindent
Since checking for a schedule with makespan~\(L\)
is equivalent to checking whether all jobs can meet a due date of~\(L\),
all \NP{}\hyp hardness results
from
\(\pJ||\Sigma C_j\) and \(\pO||\Sigma C_j\)
in \cref{sec:jcmax}
carry over to
\(\pJ||\Sigma U_j\) and \(\pO||\Sigma U_j\).
In fact,
these throughput maximization variants
turn out to be even harder:
now,
even the special cases  \(\pJ2|p_{ij}{=}1|\Sigma w_jU_j\)
and \(\pJ2|p_{ij}{=}1,r_j|\Sigma U_j\) on two machines
are NP-hard \citep{Kravchenko2000b,
  Timkovsky1998}.

The situation is more relaxed for open shop scheduling:
\(\pO|p_{ij}{=}1,r_j|\Sigma U_j\)
is \NP{}-hard for an unbounded number of machines \citep{Kravchenko2000},
whereas
\(\pO|p_{ij}{=}1|\Sigma w_jU_j\)
is even fixed\hyp parameter tractable parameterized by
the number~\(m\) of machines \citep{BruckerEtAl1993}.
For %
\(\pO m|p_{ij}{=}1,r_j|\Sigma w_jU_j\),
we only know %
polynomial\hyp time solvability
for a \emph{constant} number of machines
\citep{Baptiste2003}.

\begin{openproblem}
  Is \(\pO|p_{ij}{=}1,r_j|\Sigma w_jU_j\) fixed\hyp parameter tractable parameterized
  by the number of machines?
\end{openproblem}

\noindent
It is known that the special case \(\pO|p_{ij}{=}1,r_j|C_{\max}\)
is polynomial\hyp time solvable
and
that \(\pO|p_{ij}{=}1,r_j|\Sigma w_jU_j\)
is equivalent
to the problem \(\pP|p_{j}{=}m,\pmtn{}_{\mathbb Z},r_j|\Sigma w_jU_j\)
of scheduling jobs with processing time~\(m\)
on \(m\)~parallel identical machines
with release dates
and preemption allowed at integer times \citep{BruckerEtAl1993b}.

\paragraph{Acknowledgements}
Matthias Mnich is supported by DFG grant 
MN~59/4-1 and ERC Starting Grant 306465 (BeyondWorstCase);
he expresses his gratitude towards Alexander Grigoriev for enlightening discussions and many helpful comments.

René van Bevern is supported
by the Russian Foundation for Basic Research,
  project 16-31-60007 mol\textunderscore a\textunderscore dk,
  and by the
  Ministry of Science and Education of the Russian Federation
  under the 5-100 Excellence Program;
  he
  thanks Gerhard J.\ Woeginger for pointing out \cref{prob:o3}.

  Both authors thank Martin Kouteck{\'y} for pointing out that
  $P||C_{\max}$ is
  fixed\hyp parameter tractable for the number of distinct processing times when those are encoded in unary
  and Gerhard J.\ Woeginger for pointing out that
  \(\pO 3||C_{\max}\)~is fixed\hyp parameter tractable
  when parameterized by the maximum processing time.
They also thank the reviewers for their detailed and helpful remarks which lead to an improved presentation.

\bibliographystyle{plainnatfj} %
\bibliography{fptscheduling} %

\end{document}